\documentclass[10pt]{amsart}
\usepackage{amsmath,amssymb}
\usepackage[arrow,matrix,curve]{xy}

\usepackage{mathrsfs,bbm,eucal}
\theoremstyle{plain}
\newtheorem{thm}{Theorem}
\theoremstyle{plain}
\newtheorem{lem}{Lemma}
\theoremstyle{plain}
\newtheorem{prop}{Proposition}
\theoremstyle{plain}

\newtheorem{cor}{Corollary}
\theoremstyle{definition}
\newtheorem{defn}{Definition}
\theoremstyle{remark}

\newcommand{\dirac}{\mbox{$\mathcal{D}\!\!\!\!\!\:/\!\;$}}

\newcommand{\spinor}{\mbox{$S\!\!\!\!\!\:/\;\!$}}
\newcommand{\bundle}[1]{\CMcal{#1}}
\newcommand{\Cl}{\mathit{Cl}}
\newcommand{\R}{\mathbbm{R}}
\newcommand{\C}{\mathbbm{C}}
\newcommand{\K}{\mathbbm{K}}
\newcommand{\Q}{\mathbbm{Q}}

\newcommand{\Z}{\mathbbm{Z}}

\newcommand{\curv}[1]{\CMcal{#1}}

\begin{document}
\title{Scalar curvature on compact symmetric spaces}
\author{Mario Listing}
\address{Mathematisches Institut, Albert-Ludwigs-Universit\"at Freiburg, Eckerstra\ss e 1, 79104 Freiburg, Germany}
\email{mario.listing@math.uni-freiburg.de}
\thanks{Supported by the German Science Foundation}
\begin{abstract}
A classic result by Gromov and Lawson states that a Riemannian metric of non--negative scalar curvature on the Torus must be flat. The analogous rigidity result for the standard sphere was shown by Llarull. Later Goette and Semmelmann generalized it to locally symmetric spaces of compact type and nontrivial Euler characteristic. In this paper we improve the results by Llarull and Goette, Semmelmann. In fact we show that if $(M,g_0)$ is a locally symmetric space of compact type with $\chi (M)\neq 0$ and $g$ is a Riemannian metric on $M$ with $\mathrm{scal}_g\cdot g\geq \mathrm{scal}_0\cdot g_0$, then $g$ is a constant multiple of $g_0$. The previous results by Llarull and Goette, Semmelmann always needed the two inequalities $g\geq g_0$ and $\mathrm{scal}_g\geq \mathrm{scal}_0$ in order to conclude $g=g_0$. Moreover, if $(S^{2m},g_0)$ is the standard sphere, we improve this result further and show that any metric $g$ on $S^{2m}$ of scalar curvature $\mathrm{scal}_g\geq (2m-1)\mathrm{tr}_g(g_0)$ is a constant multiple of $g_0$.
\end{abstract}
\maketitle

\section{Introduction}

For closed surfaces the Gau\ss --Bonnet formula implies strong relations between the topology and the geometry. For instance, any metric on the $2$--torus which has non--negative or non--positive scalar curvature (Gau\ss \ curvature) is flat. Although the generalized Gau\ss --Bonnet theorem does not provide this kind of result in dimension $n\geq 3$, there are similar scalar curvature rigidity phenomena. A first example of scalar curvature rigidity in dimension $n\geq 3$ was proven by Gromov and Lawson in \cite{GrLa2}. They showed that any metric of non--negative scalar curvature on the torus has to be flat. In contrast to the $2$--dimensional case this result does not hold under the assumption of non--positive scalar curvature, since any manifold of dimension $n\geq 3$ carries a metric of constant negative scalar curvature (cf.~\cite{Loh2} for instance). In order to get scalar curvature rigidity of closed manifolds with non--vanishing Yamabe constant, one has to assume a certain scaling invariant condition on the Riemannian metric. Motivated by Gromov's $\mathrm{K-area}$ inequality in \cite{Gr01} and the fact that the $\mathrm{K-area}$ satisfies $\mathrm{K-area }(M,g)\geq \mathrm{K-area }(M,h)$ for any $g,h$ with $g\geq h$ on $\Lambda ^2TM$, Llarull used in \cite{Llarull} the algebraic condition $g\geq g_0$ to show scalar curvature rigidity of the standard sphere. In particular, he proved that if $(M,g)$ is a compact spin manifold and $f:M\to S^n$ is a smooth map of non--zero degree such that $g\geq f^*g_{S^n}$ on $\Lambda ^2TM$, then there is a point $p\in M$ with $\mathrm{scal}(g)<n(n-1)$ or $(M,g)$ is isometric to $S^n$. The analogous rigidity result for the complex and quaternionic projective space was proved by Kramer in \cite{Kra}. Min--Oo showed in \cite{MinO2} a similar result for Hermitian symmetric spaces. The general case, i.e.~the problem for symmetric spaces of compact type, was considered by Goette and Semmelmann in \cite{GoSe1}. They proved in particular, that if $(M_0,g_0)$ is a closed manifold with curvature operator $\curv{R}^0\geq 0$, $\mathrm{Ric}_0>0$ and Euler characteristic $\chi (M_0)\neq 0$, then any area non--increasing spin map $f:M\to M_0$ with  $\deg _{\widehat{A}}(f)\neq 0$ must be a Riemannian submersion if the scalar curvature satisfies
\[
\mathrm{scal}_g\geq \mathrm{scal}_0\circ f.
\]
Hence, if $(M_0,g_0)=G/H$ is a locally symmetric space of compact type with nontrivial Euler characteristic or nontrivial Kervaire semi characteristic, then
\begin{equation}
\label{gse}
g\geq g_0 \ \ \text{on}\ \ \Lambda ^2TM_0 \quad \text{and}\quad \mathrm{scal}_g\geq \mathrm{scal}_0 
\end{equation}
implies $g=g_0$. Although Goette and Semmelmann also considered the case $\mathrm{rk}(G)-\mathrm{rk}(H)=1$, there is a mistake in their proof if the Kervaire semi characteristic $\sigma (G/H)$ vanishes (cf.~\cite{Go1}). Goette generalized in \cite{Go1} the results for symmetric spaces to normal homogeneous spaces $G/H$ with $\chi (G/H)\neq 0$ respectively $\sigma (G/H)\neq 0$.

In this paper we omit the assumption of area (or distance) non--increasing maps which was always made in the previous results. We define instead a certain non--negative function on the manifold and show that the scalar curvature cannot be larger than this function. For instance as a corollary we prove that if $g$ is a Riemannian metric on a locally symmetric space $(M_0,g_0)$ of compact type with $\chi (M_0)\neq 0$ (respectively $\sigma (M_0)\neq 0$) and
\begin{equation}
\label{ineq0}
\mathrm{scal}_g\cdot g\geq \mathrm{scal}_0\cdot g_0 \qquad \text{on}\ \ TM_0,
\end{equation}
then $g$ is a constant multiple of $g_0$. In particular, comparing this inequality with the inequalities in (\ref{gse}), the assumptions in this paper are localized and much weaker. Moreover, if the model space $(M_0^n,g_0)$, $n=2m$ or $n=4k+1$, is the standard sphere of constant sectional curvature $K=1$ we can further improve inequality (\ref{ineq0}) and show that any metric $g$ on $S^n$ such that
\begin{equation}
\label{inequa0}
\mathrm{scal}_g\geq (n-1)\cdot \mathrm{tr}_g(g_0)
\end{equation}
is a multiple of $g_0$. This inequality is simply the trace of (\ref{ineq0}), i.e.~far weaker than the above condition. As we shall see, inequalities (\ref{ineq0}) and (\ref{inequa0}) have further generalizations to comparision of $2$-forms respectively the $\Lambda ^2$--trace.

A smooth map $f:M\to N$ between oriented manifolds is called a \emph{spin map} if the second Stiefel Whitney classes are related by
\[
w _2(TM)=f^*w _2(TN).
\]
Thus, in case $N$ is spin, $f$ is a spin map if and only if $M$ is spin manifold. Suppose $g$ is a Riemannian metric on $M$ and $h$ is a Riemannian metric on $N$, then $f$ is said to be \emph{distance non--increasing} at $p\in M$ if $|\mathrm{d}f(v)|_h\leq |v|_g$ for all $v\in T_pM$ and $f$ is said to be \emph{area non--increasing} at $p$ if $| \mathrm{d}f(v)\wedge \mathrm{d}f(w)| _h\leq | v\wedge w| _g$ for all $v,w\in T_pM$. Let
\[
\mathrm{dist}(f): M\to [0,\infty ), \ p\mapsto  \max _{0\neq v \in T_pM} \frac{f^*h(v , v )}{g(v ,v )}
\]
be the non--negative function on $M$ which defines the maximal distance scaling by $f$ at each point of $M$ and let
\[
\mathrm{area}(f): M\to [0,\infty ), \ p\mapsto  \max _{0\neq \eta \in \Lambda ^2T_pM} \frac{f^*h(\eta , \eta )}{g(\eta ,\eta )}
\]
be the function on $M$ which defines the maximal area scaling by $f$ at each point. In particular, $f$ is distance non--increasing for all $p\in M$ if and only if $\mathrm{dist}(f)\leq 1$ and $f$ is area non--increasing for all $p\in M$ if and only if $\mathrm{area}(f)\leq 1$. As expected, the area scaling can be estimated by the square of the length scaling:
\[
0\leq \mathrm{area}(f)\leq \mathrm{dist}(f)^2.
\]
The following theorem generalizes previous results by omitting the assumption $\mathrm{area}(f)\leq 1$ and using the function $\mathrm{area}(f)$ in the scalar curvature inequality. Note that this is a significant improvement since at points where $\mathrm{area}(f)$ vanishes, the scalar curvature inequality reads $\mathrm{scal}_g\geq 0$ while in \cite{GoSe1,Llarull} the scalar curvature is assumed to be greater than the scalar curvature on $M_0$. 
\begin{thm}
Let $(M_0^n,g_0)$, $n\geq 3$, be an oriented closed Riemannian manifold with non--negative curvature operator $\curv{R}^0\geq 0$, Ricci curvature $\mathrm{Ric}_0> 0$ and Euler characteristic $\chi (M_0)\neq 0$. Suppose $(M,g)$ is an oriented closed Riemannian manifold and $f:M\to M_0$ is a spin map of nonzero degree. If the scalar curvature satisfies
\begin{equation}
\label{ineq1}
\mathrm{scal}_g\geq  \bigl( \mathrm{scal}_0\circ f\bigl) \cdot \sqrt{\mathrm{area}(f)},
\end{equation}
then $\alpha :=\sqrt{\mathrm{area}(f)}$ is constant and $f:(M,\alpha \cdot g)\to (M_0,g_0)$ is a Riemannian covering. 

Replacing the function $\sqrt{\mathrm{area}(f)}$ by $\mathrm{dist}(f)$ the statement holds too in the case $\dim M_0=2$.
\end{thm}
In the last section we give a generalization of this theorem which assumes $\deg _{\widehat{A}}(f)\neq 0$ instead of $\deg (f)\neq 0$, in particular $\dim M=\dim M_0+4k$. However, in this case we can only show that $f:(M,g)\to (M_0,g_0)$ is a conformal submersion. In order to show that $f$ is a Riemannian submersion (like in \cite{GoSe1}) additional assumptions are necessary. If we consider the case $M=M_0$ and $f=\mathrm{id}$ in the above theorem, then inequality (\ref{ineq1}) is equivalent to
\begin{equation}
\label{ineq2}
\mathrm{scal}_g\geq \mathrm{scal}_0\cdot \| g_0\| _g
\end{equation}
where $\| g_0\| _g:M\to [0,\infty )$ is the function on $M$ defined by
\[
\| g_0\| _g(p):=\sqrt{\max _{0\neq \eta \in \Lambda ^2T_pM} \frac{g_0(\eta ,\eta )}{g(\eta ,\eta )}}\ .
\]
In fact, if a metric $g$ on $M_0$ satisfies inequality (\ref{ineq2}) or (\ref{ineq0}), then $g$ is a constant multiple of $g_0$. Note that inequality (\ref{ineq0}) implies (\ref{ineq2}), but in general the converse is not true. In order to give a $\Lambda ^2$--analog of inequality (\ref{ineq0}) we have to assume that $\mathrm{scal}_g\geq 0$, in fact (\ref{ineq2}) is equivalent to
\[
\mathrm{scal}_g\geq 0 \quad \text{and}\quad (\mathrm{scal}_g)^2\cdot g\geq (\mathrm{scal}_0)^2\cdot g_0 \quad \text{on}\ \  \Lambda ^2TM_0.
\] 
This corollary of the theorem has an odd dimensional analog. A topological version of the following theorem could also be proved, but this means to make a very complicated $KO$--theoretic assumption on the map $f:M\to M_0$. Remember that the \emph{Kervaire semicharacteristic} of a $4k+1$ dimensional manifold $M$ is given by $\sigma (M)=\sum _jb_{2j}(M)\mod 2$ where $b_i$ stands for the $i$-th Betti number of $M$. We set $\sigma (M)=0$ if $\dim M\not\equiv 1\mod 4$.
\begin{thm}
Suppose $(M,g_0)$ is an oriented Riemannian manifold of dimension $n=4k+1$ with $\curv{R}^0\geq 0$, $\mathrm{Ric}_0>0$ and $\sigma (M)\neq 0$. If $g$ is a Riemannian metric on $M$ with scalar curvature
\[
\mathrm{scal}_g\geq \mathrm{scal}_0\cdot \| g_0\| _g\ ,
\]
then $g$ is a constant multiple of $g_0$.
\end{thm}
\begin{cor}
Suppose $(M_0,g_0)$, $\dim M_0\geq 3$, is an oriented closed locally symmetric space with $\mathrm{Ric}_0>0$ and $\chi (M_0)\neq 0$ or $\sigma (M_0)\neq 0$. Let $g$ be another Riemannian metric on $M_0$. If
\[
\mathrm{scal}_g\geq \mathrm{scal}_0\cdot \| g_0 \| _g,
\]
then $g$ is a constant multiple of $g_0$.
\end{cor}
We introduced above the maximal distance and the maximal area scaling by $f$. There are of course notions of average distance respectively area scaling by $f$. We define as average distance scaling of $f:(M^n,g)\to (N,h)$ the function
\[
\mathrm{tr}_g(f^*h)=\sum _{i=1}^n(f^*h)(e_i,e_i):M\to [0,\infty )
\]
and as average area scaling
\[
\mathrm{tr}_g\left( f^*h^{\Lambda ^2}\right) =\sum _{1\leq i<j\leq n}(f^*h\otimes f^*h)(e_i\wedge e_j,e_i\wedge e_j):M\to [0,\infty ),
\]
in this case $e_1,\ldots ,e_n$ is an $g$--orthonormal basis of $TM$. Obviousely, the following inequalities hold
\[
\frac{1}{n}\mathrm{tr}_g(f^*h)\leq \mathrm{dist}(f) ,\qquad \frac{2}{n(n-1)}\mathrm{tr}_g\left( f^*h^{\Lambda ^2}\right) \leq \mathrm{area}(f)
\]
with equality (in one equation) if and only if $f^*h=g$ [in case $\frac{2}{n(n-1)}\mathrm{tr}_g\left( f^*h^{\Lambda ^2}\right) = \mathrm{area}(f)$ we need additionally $n\geq 3$ to conclude $f^*h=g$]. A simple comparison of the eigenvalues proves
\[
\mathrm{tr}_g\left( f^*h^{\Lambda ^2}\right) \leq \frac{n-1}{2n}\left[ \mathrm{tr}_g(f^*h)\right] ^2
\]
and together with the theorem below, we obtain inequality (\ref{inequa0}). The following theorem represents a significant improvement to theorem 1 and compared to Llarull's rigidity result it needs far weaker assumptions on the scalar curvature. However, it is uncertain if we can generalize this result to symmetric spaces. For the proof of theorem 3 the very simple structure of the curvature operator on $S^n$ is essential. 
\begin{thm}
Let $M$ be a closed spin manifold and $f:M\to S^{n}$, $n=2m\geq 4$, be a map of nonzero degree. If there is a metric $g$ on $M$ of scalar curvature
\begin{equation}
\label{inequal}
\mathrm{scal}_g\geq \sqrt{2n(n-1)\cdot \mathrm{tr}_g\left( f^*g_0^{\Lambda ^2}\right) },
\end{equation}
then $f:(M,g)\to (S^n,g_0)$ is a homothety, i.e.~up to scaling by a constant an isometry (here $g_0$ is of constant sectional curvature $K=1$).
\end{thm}
\section{Curvature inequality}
In order to get a simple expression for the index of the Dirac operator in terms of characteristic classes, we are using the approach presented in \cite{GoSe1}. In this section we assume that $(M_0^m,g_0)$ is an oriented manifold with non--negative curvature operator $\curv{R}^0:\Lambda ^2TM_0\to \Lambda ^2TM_0$, in particular $g_0$ has non--negative sectional curvature $K\geq 0$. Suppose $(M,g)$ is an oriented manifold of dimension $n$ and $f:M^n\to M_0^m$ is a spin map, then the vector bundle
\[
\bundle{E}:=TM\oplus f^*TM_0
\]
admits a spin structure. The Levi--Civita connection of $g$ and the Levi--Civita connection of $g_0$ induce the connection $\nabla ^\bundle{E}=\nabla \oplus f^*\nabla ^0$ on $\bundle{E}$ which is Riemannian with respect to $g\oplus f^*g_0$. The Clifford bundle of $\bundle{E}$ is given by
\[
\Cl _\K (\bundle{E})=\Cl _\K (TM)\widehat{\otimes } f^*\Cl _\K (TM_0).
\]
with $\K =\R ,\C $. Since the $\mathrm{SO}$--frame bundle of $\bundle{E}$ has structure group $\mathrm{SO}(n)\times \mathrm{SO}(m)$, the structure group of a spin structure on $\bundle{E}$ is reducible to
\[
\mathrm{Spin}(n)\cdot \mathrm{Spin}(m):=\mathrm{Spin}(n)\times \mathrm{Spin}(m)/\{ \pm 1\} \subset \mathrm{Spin}(n+m).
\]
Let $\spinor \bundle{E}$ be the a complex (or real) spinor bundle of $\bundle{E}$ induced by a choice of the spin structure and the tensor product of the complex (respectively real) spin representations. The connection $\nabla ^\bundle{E}$ lifts uniquely to a Riemannian connection $\overline{\nabla }$ on $\spinor \bundle{E}$. For each point $p\in M$ there are neighborhoods $p\in U\subset M$ and $f(p)\in V\subset M_0$ in such a way that the spinor bundle decomposes as
\begin{equation}
\label{zerleg}
\spinor \bundle{E}_{|U}=\spinor U\otimes f^*(\spinor V)
\end{equation}
(note that we do not assume $\spinor \bundle{E}$ to be irreducible at this point). In particular, if $M_0$ is spin, then $M$ is spin and we conclude $\spinor \bundle{E}=\spinor M\widehat{\otimes }f^*\spinor M_0$ as $\Cl _\K (\bundle{E})$ module. Since the forthcoming computations are of local nature we use the notation $\spinor \bundle{E}=\spinor M\otimes f^*\spinor M_0$ even if $M_0$ is not spin. Because $\spinor \bundle{E}$ is a $\Cl _\K (\bundle{E})$--module, $\spinor \bundle{E}$ becomes a Dirac bundle over $M$ if we use the imbedding
\[
\Cl _\K (TM)\hookrightarrow \Cl _\K (TM)\otimes \mathbbm{1}\subset \Cl _\K (\bundle{E}).
\]
The corresponding Dirac operator will be denoted by
\[
\overline{\dirac }:=\sum _{j=1}^n\gamma (e_j)\overline{\nabla }_{e_j}
\]
where $e_1,\ldots ,e_n$ is an orthonormal basis of $T_pM$. Since the connection $\overline{\nabla }$ preservers the decomposition in (\ref{zerleg}), i.e.~$\overline{\nabla }=\nabla \otimes \mathbbm{1}+\mathbbm{1}\otimes f^*\nabla ^0$ on the formal tensor product $\spinor \bundle{E}=\spinor M\otimes f^*\spinor M_0$, we can use \cite[Ch.~II Thm.~8.17]{LaMi} to compute the Bochner--Weitzenb\"ock formula of $\overline{\dirac }$:
\begin{equation}
\label{Lichner}
\overline{\dirac }^2 =\overline{\nabla }^* \overline{\nabla }+\frac{\mathrm{scal}_g}{4}+\frak{R}
\end{equation}
where $\frak{R}$ is defined by
\[
\frak{R}=\frac{1}{2}\sum _{i,j=1}^n\gamma (e_i)\gamma (e_j)\otimes R_{e_i,e_j}^0
\]
and $R^0$ means the curvature of $f^*\nabla ^0$.
\begin{defn}
Suppose $(V,\langle .,.\rangle _V)$ and $(W,\langle .,.\rangle _W)$ are inner product spaces and $\beta :V\to W$ is a linear transformation. Then $\beta $ is called a \emph{homothetic injection} if there is a constant $c>0$ such that
\begin{equation}
\label{homothetic}
\left< v,\hat v\right> _V=c \cdot \left< \beta (v),\beta (\hat v) \right> _W 
\end{equation}
for all $v,\hat v\in V$. Let $V=\ker (\beta )\oplus V^\prime $ be the orthogonal decomposition of $V$ w.r.t.~the inner product. Then $\beta $ is said to be a \emph{homothetic surjection} if (\ref{homothetic}) holds for all $v,\hat v\in V^\prime $. A homothetic isomorphism is also called a \emph{homothety}.
\end{defn}
\begin{prop} 
\label{main_prop}
Let the Riemannian curvature operator of $(M_0,g_0)$ be non--negative and $\dim M_0\geq 3$, then the curvature endomorphism $\frak{R}$ is at each point bounded as follows
\begin{equation}
\label{ineq}
\frak{R}\geq -\frac{\mathrm{scal}_0\circ f}{4}\sqrt{\mathrm{area}(f)}
\end{equation}
where $\mathrm{scal}_0$ means the scalar curvature of $g_0$. Furthermore, suppose $\mathrm{Ric}(g_0)$ is positive definite at $f(p)\in M_0$ and $-\frac{1}{4}\sqrt{\mathrm{area}(f)} \cdot \mathrm{scal}_0( f(p))$ is the minimal eigenvalue of $\frak{R}$ at $p\in M$, then $f_*:T_pM\to T_{f(p)}M_0$ is a homothetic surjection or $\mathrm{area}(f)(p)=0$. 

In particular, if $\mathrm{Ric}(g_0)>0$ on $M_0$ and $U\subset M$ denotes the interior of all points  $p\in M$ where the minimal eigenvalue of $\frak{R}$ is $-\frac{1}{4}\sqrt{\mathrm{area}(f)}\cdot \mathrm{scal}_0(f(p))$ and $\mathrm{area}(f)(p)>0 $, then
\[
f:(U,\sqrt{\mathrm{area}(f)}\cdot g) \to (f(U),g_0)
\]
is a Riemannian submersion. 
\end{prop}
Note that $\sqrt{\mathrm{area}(f)}$ is smooth on $U$, while of course $U$ could be empty in this proposition. Recall that if $\mathrm{area}(f)$ vanishes at $p$, then the image of $f_*:T_pM\to T_{f(p)}M_0$ is at most one dimensional and $\mathrm{Im}(f_*)$ is trivial if and only if $\mathrm{dist}(f)=0$. The assumption $\dim M_0\geq 3$ can be omitted in the above proposition if we replace the function $\sqrt{\mathrm{area}(f)}$ everywhere by $\mathrm{dist}(f)$ (in case $\dim M_0=2$, $\Lambda ^2M_0$ has rank one which will not be enough to show that $f_*$ is a homothetic surjection if $\mathrm{area}(f)>0$). The following two example provide non--constant maps which are not homothetic surjections but satisfy all the assumptions of the proposition except $\mathrm{Ric}(g_0)>0$ respectively $\mathrm{area}(f)>0$ ($T^n$ means the standard $n$--dimensional flat torus, $c:[0,2\pi )\to S^n$ is a simple closed geodesic):
\[
f:T^n\to S^n,\ (t_1,\ldots ,t_n)\mapsto c(t_1)\ , \qquad \widetilde{f}:T^n\to T^n\times S^1, \ p\mapsto (p,t_0).
\]
In both cases $\frak{R}$ vanishes and inequality (\ref{ineq}) is an equality. Although $\mathrm{Ric}(g_{S^n})$ is positive definite, $f$ is nowhere a homothetic surjection since $\mathrm{area}(f)=0$. Moreover, $\widetilde{f}$ satisfies $\mathrm{area}(f)\equiv 1$, but $\widetilde{f}$ is nowhere a homothetic surjection since $\mathrm{Ric}(g_0)=0$. In fact, in both cases we have $U=\emptyset $.
\begin{prop}
Suppose $(M_0,g_0)$ is the standard sphere $S^n$ of constant sectional curvature $K=1$ and $\dim M=n$. Then the curvature endomorphism satisfies
\begin{equation}
\label{inequal1}
\frak{R}\geq -\frac{1}{2}\sqrt{\frac{n(n-1)}{2}\mathrm{tr}_g\left( f^*g_0^{\Lambda ^2}\right) }
\end{equation}
and $\frak{R}$ attains this minimal eigenvalue at a point $p\in M$ if and only if $f_*:T_pM\to T_{f(p)}S^n$ is a homothety or $\mathrm{area}(f)(p)=0$. In particular, if $U\subset M$ is the interior of all points where the minimal eigenvalue of $\frak{R}$ is $-\frac{1}{2}\sqrt{\frac{n(n-1)}{2}\mathrm{tr}_g\left( f^*g_0^{\Lambda ^2}\right) }$ and $\mathrm{area}(f)>0$, then
\[
f:\Bigl( U, \sqrt{\mathrm{area}(f)}\cdot g\Bigl) \to (f(U),g_0)
\]
is a Riemannian covering. 
\end{prop}
To show these two propositions we will simplify the curvature expression $\frak{R}$. For each point $p\in M$ the map $f:M\to M_0$ induces an isometric isomorphism $\beta _p:f^*\spinor \frak{p}\to \spinor \frak{p}$ where $\frak{p}=T_{f(p)}M_0$ and $\spinor \bundle{E}_{p}=\spinor T_pM\otimes f^*\spinor \frak{p}$. Since the curvature of the connection $\nabla ^0$ is the curvature of the (virtual) spin bundle $\spinor M_0$, we obtain the curvature of $f^*\nabla ^0$
\begin{equation}
\label{curv12}
\begin{split}
R^0_{v,w}&=\beta ^{-1}_p\circ R^{0}_{f_*v,f_*w}\circ \beta _p\\
&=-\frac{1}{2}\beta ^{-1}_p\circ \gamma ^0(\curv{R}^0(f_*(v\wedge w)))\circ \beta _p
\end{split}
\end{equation}
where $v,w\in T_pM$, $\gamma ^0 $ is the Clifford multiplication on $\spinor \frak{p}$ and $\curv{R}^0$ is the Riemannian curvature operator of $(M_0,g_0)$ considered as endomorphism on $\Lambda ^2M_0$. Thus, the curvature operator $\frak{R}$ is determined by
\begin{equation}
\label{curv11}
\begin{split}
\frak{R}&=\frac{1}{2}\sum _{ i,j=1}^{n}\gamma (e_i)\gamma (e_j)\otimes R_{e_i,e_j}^0\\
&=-\frac{1}{4}\sum _{i,j=1}^{n}\gamma (e_i\wedge e_j)\otimes \beta _p^{-1}\gamma ^0\Bigl( \curv{R}^0(f_*e_i\wedge f_*e_j)\Bigl) \beta _p.
\end{split}
\end{equation}
Let $B\in \Gamma (\mathrm{End}(TM))$ be the symmetric positive semi--definite transformation defined by
\begin{equation}
\label{defn_of_B}
\begin{split}
g(BX,Y)&=f^*g_0(X,Y)=g_0(f_*(X),f_*(Y))
%g(BX,Y)&=g(X,BY),
\end{split}
\end{equation}
and set $B_k:=\underbrace{B\otimes \cdots \otimes B}_{k-\text{times}}\in \Gamma (\mathrm{End}(\Lambda ^kTM))$, then $B=B_1$ and $B_2$ satisfy
\[
0\leq B\leq \mathrm{dist}(f)\ ,\qquad 0\leq B_2\leq \mathrm{area}(f)\leq \mathrm{dist}(f)^2
\]
(note that the upper inequalities are sharp in each point, since by definition, $B$ has an eigenvalue $\mathrm{dist}(f)$ and $B_2$ has an eigenvalue $\mathrm{area}(f)$). Moreover, we obtain by definition and a standard computation
\[
\mathrm{tr}(B)=\mathrm{tr}_g(f^*g_0),\quad \mathrm{tr}(B_2)=\mathrm{tr}_g\left( f^*g_0^{\Lambda ^2}\right) , \quad [\mathrm{tr}(B)]^2=2\cdot \mathrm{tr}(B_2)+\mathrm{tr}(B^2).
\]
In fact, if $\dim M_0=\dim M=n$, the Cauchy Schwarz inequality yields $[\mathrm{tr}(B)]^2\leq n\cdot \mathrm{tr}(B^2)$ with equality iff $B$ is a multiple of $\mathrm{Id}$, i.e.~we conclude
\begin{equation}
\label{trace_formula}
\mathrm{tr}(B_2)\leq \frac{n-1}{2n}[\mathrm{tr}(B)]^2
\end{equation}
with equality if and only if $B$ is a multiple of $\mathrm{Id}$. Note altough $\mathrm{tr}(B)\leq n\cdot \mathrm{dist}(f)$, there is no estimate of $\mathrm{tr}(B)$ in terms of $\mathrm{area}(f)$ (if $B$ has only one nonvanishing eigenvalue, $\mathrm{area}(f)=0$ while $\mathrm{tr}(B)>0$). Conversely, there is no sharp (respectively optimal) estimate of $\mathrm{area}(f)$ in terms of $\mathrm{tr}(B)$.
\begin{lem}
\label{lem_lin}
Suppose $f:M\to N$ is a differentiable map, $g$ is a Riemannian metric on $M$ and $\overline{g}$ is a Riemannian metric on $N$. Then $g$ and $\overline{g}$ induce isomorphisms $t:\Lambda ^kT^*_pM\rightarrow \Lambda ^kT_pM$, $r:\Lambda ^kT_qN\rightarrow \Lambda ^kT_q^*N$ and the following diagram is commutative:
\[
\begin{xy}
\xymatrix{
\Lambda ^kT_pM \ar[r]^{f_*} \ar[rd]^{B_k}& \Lambda ^kT_{f(p)}N\ar@{.>}[d]^{f^\#}\ar[r]^{r} & \Lambda ^kT_{f(p)}^*N\ar[d]^{f^*}\\
& \Lambda ^kT_pM &\ar[l]^t  \Lambda ^kT_p^*M
}
\end{xy}
\]
In particular, $f^\#:\Lambda ^kT_{f(p)}N\to \Lambda ^kT_pM$ is uniquely determined for each $p\in M$ by $f^\#:=t\circ f^*\circ r$ and satisfies
$f^\#\circ f_*=B_k$ as well as
\[
\overline{g}(f_*v,w)=g(v,f^\#w)
\]
for all $v\in \Lambda ^kT_pM$ and $w\in \Lambda ^kT_{f(p)}N$. Furthermore, define
\[
\breve{B}_k:=f_*f^\# :\Lambda ^kT_{f(p)}N\to \Lambda ^kT_{f(p)}N,
\]
then $\breve{B}_k$ is positive semi--definite and symmetric w.r.t.~$\overline{g}$. Moreover, the non--vanishing eigenvalues of $B_k$ at $p\in M$ coincide with the non--vanishing eigenvalues of $\breve{B}_k$ at $f(p)\in N$ and we have $\breve{B}_k=\underbrace{\breve{B}\otimes \cdots \otimes \breve{B}}_k$ with $\breve{B}=\breve{B}_1$. 
\end{lem}
\begin{proof}
This lemma collects some facts from linear algebra. Note that $\breve{B}$ can be defined analogous to (\ref{defn_of_B}) by
\[
\begin{split}
\overline{g}(\breve{B}x,y)&=g(f^\# x,f^\# y)
%\overline{g}(\breve{B}x,y)&=\overline{g}(x,\breve{B}y),
\end{split}
\]
for all $x,y\in T_{f(p)}N$. If $v\in T_pM$ is an eigenvector of $B$ to the eigenvalue $\lambda \neq 0$, then (\ref{defn_of_B}) yields that $f_*v\neq 0$. In particular, $f_*v\in T_{f(p)}N$ is an eigenvector of $\breve{B}$ to the eigenvalue $\lambda $:
\[
\breve{B}f_*v=(f_*f^\# )f_*v=f_*(f^\# f_*)v=f_*Bv=\lambda f_*v
\]
and thus, $\lambda $ is an eigenvalue of $\breve{B}$ (appears with the same multiplicity as in $B$). That any nonzero eigenvalue of $\breve{B}$ is an eigenvalue of $B$ follows in the same way.
\end{proof}

Let $h_1,\ldots ,h_s\in \Lambda ^2T_{f(p)}M_0$ be a $g_0$--orthonormal eigenbasis of $\curv{R}^0$ and $\kappa _1,\ldots ,\kappa _s$ be the corresponding eigenvalues (note that $\curv{R}^0$ is symmetric), then $\curv{R}^0\geq 0$ yields
\[
\kappa _j=\left< \curv{R}^0(h_j),h_j\right> \geq 0.
\]
Furthermore, we obtain from (\ref{curv11}) and the symmetry of $\curv{R}^0$ ($e_1,\ldots ,e_n$ is a $g$--orthonormal basis of $T_pM$):
\begin{equation}
\label{curv13}
\begin{split}
\frak{R}&=-\frac{1}{4}\sum _{i,j=1}^n\sum _{l=1}^sg_0(f_*(e_i\wedge e_j),\curv{R}^0(h_l))\cdot \gamma (e_i\wedge e_j)\otimes \beta _p^{-1}\gamma ^0(h_l) \beta _p\\
&=-\frac{1}{2}\sum _{l=1}^s\kappa _l\gamma (f^\# h_l)\otimes \beta _p^{-1}\gamma ^0(h_l)\beta _p.
\end{split}
\end{equation}
Let $\alpha \in [0,\infty )$ be non--negative and define
\begin{equation}
\label{defnC}
\frak{C}:=\sum _{l=1}^s\kappa _l[\gamma (f^\# h_l)\otimes \mathrm{Id}+\alpha \cdot \mathrm{Id}\otimes \beta _p^{-1}\gamma ^0(h_l)\beta _p]^2 \in \mathrm{End}(\spinor \bundle{E}_p).
\end{equation}
then a straightforward calculation shows
\[
\frak{C}=\sum _{l=1}^s\kappa _l[(\gamma (f^\# h_l))^2\otimes \mathrm{Id}+\alpha ^2\cdot \mathrm{Id}\otimes \beta _p^{-1}(\gamma ^0(h_l))^2\beta _p]-4\alpha \frak{R}.
\]
\begin{lem}
Suppose $\eta $ is a $2$--form, then
\[
\gamma (\eta )^2=\gamma (\eta \wedge \eta )-|\eta |^2
\]
for any Clifford module. Moreover, let $\kappa _1,\ldots ,\kappa _s$ be the eigenvalues of $\curv{R}^0$ and $h_1,\ldots ,h_s$ be the corresponding orthonormal eigenbasis, then $\sum \kappa _l=\frac{1}{2}\mathrm{scal}_0$ and
\[
\sum _{l=1}^s\kappa _lh_l\wedge h_l=0.
\]
\end{lem}
\begin{proof}
The first statement is a straightforward calculation and the second follows from $\mathrm{scal}_0=2\mathrm{trace}(\curv{R}^0)$. Moreover, the first Bianchi identity yields after applying vectors $x,y,z,t$:
\[
\sum _{l=1}^s\kappa _lh_l\wedge h_l=\sum _l\curv{R}^0(h_l)\wedge h_l=0.
\]
\end{proof}
\begin{proof}[Proof of Proposition 1]
Using the fact $f^\# (h_l\wedge h_l)=f^\# h_l\wedge f^\# h_l$, this lemma simplifies $\frak{C}$ to:
\begin{equation}
\label{simplC}
\frak{C}=-4\alpha \frak{R}-\frac{\alpha ^2}{2}\cdot \mathrm{scal}_0(f(p))-\sum _{l=1}^s\kappa _l |f^\#h_l|^2_g.
\end{equation}
Since $0\leq B_2\leq \mathrm{area}(f)$ and the nonzero eigenvalues of $B_2$ and $\breve{B}_2$ coincide, lemma \ref{lem_lin} yields
\[
g(f^\# h_l,f^\# h_l)=g_0(f_*f^\#h_l,h_l)=g_0(\breve{B}_2h_l,h_l)\leq \mathrm{area}(f).
\]
In particular, $\kappa _j\geq 0$ and $\sum _l\kappa _l=\frac{1}{2}\mathrm{scal}_0(f(p))$ prove the following inequality
\begin{equation}
\label{ineq_loc}
\sum _{l=1}^s\kappa _l|f^\# h_l|^2_g \leq \frac{1}{2}\mathrm{area}(f)\cdot \mathrm{scal}_0(f(p)).
\end{equation}
This estimate completes the proof for the first part of proposition \ref{main_prop}. Consider the definition of $\frak{C}$ in (\ref{defnC}). Since $\gamma (\eta )$ is a skew adjoint action on any Clifford module for arbitrary $2$--forms $\eta $ and $\kappa _l\geq 0$ for all $l=1\ldots s$, we conclude from (\ref{defnC}): $\frak{C}\leq 0$. Hence, set $\alpha :=\sqrt{\mathrm{area}(f)}:M\to [0,\infty )$ then equation (\ref{simplC}) and inequality (\ref{ineq_loc}) show
\begin{equation}
\label{ineq_fin}
\frak{R}\geq -\frac{\alpha }{8}\mathrm{scal}_0\circ f-\frac{1}{4\alpha }\sum _{l=1}^m\kappa _l|f^\# h_l|^2_g\geq -\frac{\alpha }{4}\mathrm{scal}_0\circ f
\end{equation}
(note that in case $\alpha =\sqrt{\mathrm{area}(f)}=0$ at a point $p\in M$, $f_*:\Lambda ^2T_pM\to \Lambda ^2T_{f(p)}M_0$ has to vanish and thus, $\frak{R}$ vanishes at $p\in M$ from (\ref{curv11}), i.e.~this inequality is also true at points where $\alpha = 0$). In order to show the second part of proposition \ref{main_prop} we need the conditions $\mathrm{Ric}(g_0)>0$ and $m=\dim M_0\geq 3$. We have to prove that $f_*:T_pM\to T_{f(p)}M_0$ is a homothetic surjection if $\mathrm{area}(f)(p)>0$. Suppose the minimal eigenvalue of $\frak{R}$ at $p\in M$ is $-\frac{1}{4}\alpha \cdot \mathrm{scal}_0(f(p))$ with $\alpha :=\sqrt{\mathrm{area}(f)(p)}>0$.  In this case we obtain equality for at least one nontrivial spinor in (\ref{ineq_fin}) and henceforth, we obtain equality in (\ref{ineq_loc}). In particular, $| f^\# h_l| _g=\alpha $ for all $l$ with $\kappa _l> 0$ which is equivalent to $\breve {B}_2 =\alpha ^2\mathrm{Id}$ on $\mathrm{Im}(\curv{R}^0)\subset \Lambda ^2T_{f(p)}M_0$. Let $e_1,\ldots ,e_m\in T_{f(p)}M_0$ be an orthonormal eigenbasis of $\breve{B}=\breve{B}_1$ to eigenvalues $\lambda _1\geq \ldots \geq \lambda _m$. Then $\lambda _i\lambda _j$ ($i\neq j$) are the eigenvalues of $\breve{B}_2$ and $0\leq \breve{B}_2\leq \alpha ^2$ yields $\lambda _i\lambda _j\leq \alpha ^2 $ for all $i,j=1\ldots m$ with $i\neq j$. Moreover, $\breve{B}_2=\alpha ^2$ on $\mathrm{Im}(\curv{R}^0)$ supplies for the curvature of $(M_0,g_0)$ (use the symmetry of $\breve{B}_2=\breve{B}\otimes \breve{B}$):
\[
\begin{split}
\alpha ^2\cdot R^0(e_i,e_j,e_k,e_l)&=\alpha ^2\cdot g_0(\curv{R}^0(e_i\wedge e_j ),e_l \wedge e_k )=g_0(\breve{B}_2\curv{R}^0(e_i\wedge e_j),e_l\wedge e_k)\\
&=g_0(\curv{R}^0(e_i\wedge e_j),\breve{B}e_l\wedge \breve{B}e_k)=\lambda _k\lambda _lR^0(e_i,e_j,e_k,e_l).
\end{split}
\]
Since $\mathrm{Ric}(g_0)>0$ at $f(p)$, for any $k$ there is some $l$ with $R^0(e_l,e_k,e_k,e_l)>0$. Thus, for all $k=1\ldots m$ there is some $l\neq k$ with $\lambda _k\lambda _l=\alpha ^2$. We assumed $\alpha >0$ and $m=\dim M_0\geq 3$. Thus, let $k$ be arbitrary and $ l\neq k$ in such a way that $\lambda _k\lambda _l=\alpha ^2$. Suppose $i\neq k$ as well as $i\neq l$. Then
\[
\lambda _i\lambda _k\leq \alpha ^2,\qquad \lambda _i\lambda _l\leq \alpha ^2
\]
together with $\lambda _k\lambda _l=\alpha ^2$ yields $\lambda _i\leq \alpha $ and since $k$ was arbitrary, we conclude $\lambda _i\leq \alpha $ for all $i=1\ldots m$. Because for any $i$ there is some $j$ with $\lambda _i\lambda _j=\alpha ^2$, we obtain $\lambda _i=\alpha $ for all $i=1\ldots m$ which is equivalent to $\breve{B}=\alpha \mathrm{Id}$. Since the non--vanishing eigenvalues of $B$ and $\breve{B}$ coincide (lemma \ref{lem_lin}), definition \ref{defn_of_B} proves that $f_*:T_pM\to T_{f(p)}M_0$ is a homothetic surjection in case $\mathrm{area}(f)(p)>0$ (we have $f_*$ is surjective, $f^\# $ is injective, $T_pM=\ker (f_*)\oplus V$ and $B=\alpha $ on $V=\mathrm{Im}(f^\# )$ as well as $B=0$ on $\ker (f_*)$). 

Suppose now that $\mathrm{Ric}(g_0)>0$. Define $U\subset M$ to be the interior of all point $p\in M$ where the minimal eigenvalue of $\frak{R}$ is $-\frac{1}{4}\alpha \cdot \mathrm{scal}_0(f(p))$ and where $\mathrm{area}(f)(p)>0$. Then the above considerations show that $f:U\to f(U)$ is a submersion and definition \ref{defn_of_B} ($B=0$ on $\ker (f_*) \subset TU$ as well as $B=\alpha $ on $\ker (f_*)^\perp \subset TU$) prove that 
\[
f:(U,\alpha g)\to (f(U),g_0)
\]
is a Riemannian submersion with $\alpha =\sqrt{\mathrm{area}(f)}$.
\end{proof}
\begin{proof}[Proof of Proposition 2] Consider the expression of $\frak{C}$ in equation (\ref{simplC}). Since $(M_0,g_0)=S^n$ is of constant sectional curvature $1$, we know $\kappa _j=1$ for all $j$ and
\[
\frak{C}=-4\alpha \frak{R}-\frac{\alpha ^2}{2}n(n-1)-\mathrm{tr}(B_2)
\]
[use lemma \ref{lem_lin} and $\mathrm{tr}(B_2)(p)=\mathrm{tr}(\breve{B}_2)(f(p))$]. Set $\alpha :=\sqrt{\frac{2}{n(n-1)}\mathrm{tr}(B_2)}$, then $\alpha $ vanishes at $p\in M$ only if $\mathrm{area}(f)(p)=0$. In this case $f_*:\Lambda ^2T_pM\to \Lambda ^2T_{f(p)}S^n$ is trivial and $\frak{R}(p)=0$ which means that the inequality in proposition 2 is fulfilled. If $\alpha >0 $, we obtain from $\frak{C}\leq 0$ the estimate
\[
\frak{R}\geq -\frac{\alpha }{8}n(n-1)-\frac{1}{4\alpha }\mathrm{tr}(B_2)=-\frac{1}{2}\sqrt{\frac{n(n-1)}{2}\mathrm{tr}(B_2)}
\]
Suppose now that $\alpha >0$ and that the minimal eigenvalue of $\frak {R}$ at $p\in M$ is determined by $-\frac{1}{2}\sqrt{\frac{n(n-1)}{2}\mathrm{tr}(B_2)}$. 
In this case, the maximal eigenvalue of $\frak{C}$ in $p$ is $0$ and there is a nontrivial spinor $\varphi \in \spinor \bundle{E}_p$ with $\frak{C}\varphi =0$. Consider the definition of $\frak{C}=\sum _l C_l^2$ in equation (\ref{defnC}) where
\[
C_l=\gamma (f^{\#}h_l)\otimes \mathrm{Id}+\alpha \cdot \mathrm{Id}\otimes \beta _p^{-1}\gamma (h_l)\beta _p.
\]
Since $C_l^2\leq 0$ and $C_l$ is skew symmetric, we conclude that $C_l\varphi =0$ for all $l$. Defining
\[
D_l=\gamma (f^{\#}h_l)\otimes \mathrm{Id}-\alpha \cdot \mathrm{Id}\otimes \beta _p^{-1}\gamma (h_l)\beta _p
\]
yields
\[
D_lC_l= \gamma (f^{\#}h_l)^2\otimes \mathrm{Id}-\alpha ^2\mathrm{Id}\otimes \beta _p^{-1}\gamma (h_l)^2\beta _p
\]
Since $S^n$ has pure curvature operator, we can choose $h_l$ to be a wedge product $h_l=e_i\wedge e_j$ where $e_1,\ldots ,e_n\in T_{f(p)}S^n$ is an orthonormal eigenbasis of $\breve{B}$ to eigenvalues $\lambda _1,\ldots ,\lambda _n$. Thus, we obtain
\[
D_lC_l=-g_0( \breve{B}(e_i)\wedge \breve{B}(e_j),e_i\wedge e_j) +\alpha ^2=-\lambda _i\lambda _j+\alpha ^2
\]
and $D_lC_l\varphi =0$ for all $l$ shows $\lambda _i\lambda _j=\alpha ^2$ for all $i,j=1\ldots n$ and $i\neq j$. However, this proves $\lambda _i=\alpha $ for all $i=1,\ldots ,n$ if $n\geq 3$. Hence, $\breve{B}=\alpha \mathrm{Id}$ at $f(p)$ as well as $B=\alpha \mathrm{Id}$ at $p$ which proves that $f_*:T_pM\to T_{f(p)}S^n$ is a homothety. Note that $\alpha ^2=\frac{2}{n(n-1)}\mathrm{tr}(B_2)$ vanishes at $p\in M$ iff $\mathrm{area}(f)(p)=0$ and that $B=\alpha \cdot  \mathrm{Id}$ means
\[
\alpha ^2=\frac{2}{n(n-1)}\mathrm{tr}(B_2) =\mathrm{area}(f)=\mathrm{dist}(f)^2=\left[ \frac{1}{n}\mathrm{tr}(B)\right] ^2
\]
i.e.~the scaling factor is determined by $\sqrt{\mathrm{area}(f)}$.
\end{proof}
\section{Proof of the theorems in the introduction}
First we show that the kernel of the Dirac operator $\overline{\dirac }$ from the previous section is non--trivial. Consider the situation of theorem 1 and 3. Then $\deg (f)\neq 0$ and $\chi (M_0)\neq 0$ imply $\dim M=\dim M_0=2k$. Let $\spinor \bundle{E}$ be the irreducible complex spinor bundle of $\bundle{E}=TM\oplus f^*TM_0$, then $\spinor \bundle{E}$ is a complex Dirac bundle over $M$ and naturally $\Z _2$--graded by the volume form of $\Cl _\C (\bundle{E})$. The index of the Dirac operator $\overline{\dirac }^+:\Gamma (\spinor ^+\bundle{E})\to \Gamma (\spinor ^-\bundle{E})$ is given by (cf.~section 4 or proof of theorem 2.4 in \cite{GoSe1})
\[
\mathrm{ind}(\overline{\dirac }^+) =\ker \overline{\dirac }^+-\mathrm{coker} \ \overline{\dirac }^+=\chi (M_0)\cdot \deg (f).
\]
Hence, under the assumptions of theorem 1 respectively theorem 3, the kernel of $\overline{\dirac }$ is not trivial.

Now consider the situation of theorem 2. In this case we have $M=M_0$ and $f=\mathrm{id}$ in section 2. Let's start with the bundle $\bundle{E}_0=TM\oplus TM$ equipped with the product metric $g_0\oplus g_0$ and Levi--Civita connection $\nabla ^0\oplus \nabla ^0$. Then we can choose representations (depending on $\dim M=8l+1$ respectively $\dim M=8l+5$) in such a way that the real spinor space $\spinor \bundle{E}_0=\spinor M\otimes \spinor M$ is equivalent to the $\Z _2$--graded Clifford bundle $\Cl _\R (M,g_0)$. Clifford multiplication on $\spinor \bundle{E}_0$ with elements $\Cl _\R (\bundle{E}_0)=\Cl _\R (M)\widehat{\otimes }\Cl _\R (M)$ is determined on $\Cl _\R (M)$ by multiplication from the left with elements $\eta \otimes \mathbbm{1}$ and multiplication from the right with elements $\mathbbm{1}\otimes \eta $. In this way $\spinor \bundle{E}_0$ becomes a $\Z _2$--graded module for $\Cl _\R (\bundle{E}_0)$. Let $\omega _0\in \Cl _\R (M,g_0)$ be the volume form for $g_0$, then Clifford multiplication with $\Omega :=\mathbbm{1}\otimes \omega _0\in \Cl _\R ^1(\bundle{E}_0)$ on $\spinor \bundle{E}_0$ determines a parallel $\Z _2$--graded $\Cl _1$--action on $\spinor \bundle{E}_0$ (use the fact $\Omega ^2=-\mathrm{Id}$). In fact this action corresponds to right multiplication on $\Cl _\R (M,g)$ by the volume form $\omega _0$. Thus, since $M$ has dimension $4k+1$, the $\Cl _1$--index of the Dirac operator $\overline{\dirac _0}^{\mathrm{ev}}:\Gamma (\Lambda ^{\mathrm{even}})\to \Gamma (\Lambda ^{\mathrm{odd}})$ is given by (cf.~\cite[Ch.~II Example 7.7]{LaMi})
\[
\mathrm{ind}_1(\overline{\dirac _0}^\mathrm{ev})=\sigma (M).
\]
Consider the same bundle $\bundle{E}=\bundle{E}_0=TM\oplus TM$ but equipped with the metric $g\oplus g_0$ and connection $\nabla ^\bundle{E}=\nabla \oplus \nabla ^0$. Since $\Cl _\R (\bundle{E})=\Cl _\R (TM,g)\widehat{\otimes }\Cl _\R (TM,g_0)$ the form $\Omega :=\mathbbm{1}\otimes \omega _0\in \Cl _\R ^1(\bundle{E})$ is still parallel w.r.t.~the connection $\nabla ^\bundle{E}$ and hence, $\Omega $ defines again a parallel $\Z _2$--graded $\Cl _1$--action on $\spinor \bundle{E}$. In particular, the bundles $\spinor \bundle{E}$ and $\spinor \bundle{E}_0$ are equivalent $\Z _2$--graded $\Cl _1$--Dirac bundles. Thus, the $\Cl _1$--index of the Dirac operators $\overline{\dirac }^\mathrm{ev}:\Gamma (\spinor ^+\bundle{E})\to \Gamma (\spinor ^-\bundle{E})$ and $\overline{\dirac _0}^\mathrm{ev}:\Gamma (\spinor ^+\bundle{E}_0)\to \Gamma (\spinor ^-\bundle{E}_0)$ coincide and we conclude from $\sigma (M)\neq 0$ that the kernel of $\overline{\dirac }$ cannot be trivial.

Now we are able to finish the proof of the main theorems (in the following we have to set $f=\mathrm{id}$, $M=M_0$ in the case of theorem 2 and use the fact $\mathrm{scal}_0=n(n-1)$ in case of theorem 3). We define the function $\alpha :M\to [0,\infty )$ by
\renewcommand{\arraystretch}{1.5}
\[
\alpha :=\left\{ \begin{array}{lcl}\sqrt{\mathrm{area}(f)}& \ & \text{theorem 1}\\
\| g_0\| _g&\ & \text{theorem 2}\\
\sqrt{\frac{2}{n(n-1)} \mathrm{tr}_g\left( f^*g_0^{\Lambda ^2}\right) }&&\text{theorem 3}
\end{array} \right\}
\]
\renewcommand{\arraystretch}{1.0}
Consider the integrated version of (\ref{Lichner}) with $\psi \in \Gamma (\spinor \bundle{E})$:
\begin{equation}
\label{int_lich}
\int \limits _M|\overline{\dirac }\psi |^2=\int\limits _M|\overline{\nabla} \psi |^2+\frac{\mathrm{scal}_g}{4}|\psi |^2+\left< \frak{R}\psi ,\psi \right> .
\end{equation}
The inequalities (\ref{ineq}), (\ref{ineq1}) and (\ref{inequal}), (\ref{inequal1}) show that the function on the right hand side is point wise non--negative. In particular, let $0\neq \psi \in \ker (\overline{\dirac})$ be non--trivial, then $\overline{\nabla }\psi =0$ and $|\psi |^2=\mathrm{const}>0$ ($\overline{\nabla }$ is a Hermitian connection). Thus, $\psi (p)\neq 0$ for all $p\in M$ and $\left< \frak{R}\psi ,\psi \right> =-\frac{1}{4}\mathrm{scal}_g\cdot |\psi |^2$ yield
\begin{equation}
\label{scal_equal}
\mathrm{scal}_g=\alpha \cdot \mathrm{scal}_0\circ f
\end{equation}
and $\psi $ is at each point of $M$ an eigenvector of $\frak{R}$ to the eigenvalue
\[
-\frac{1}{4}\alpha \cdot \mathrm{scal}_0\circ f.
\] 
Hence, proposition \ref{main_prop} respectively proposition 2 proves that at $p\in M$ the map $f_*:T_pM\to T_{f(p)}M_0$ is a homothetic surjection or $\alpha (p)=0$. However, since $\dim M=\dim M_0$, a homothetic surjection is nothing but a homothety. Thus, let $U$ be the open set in $M$ where $\alpha >0$, then $\deg (f)\neq 0$ proves that $U$ is not empty (note that the case $\| g_0\| _g(p)=0$ cannot happen in case of theorem 2, i.e.~we have $U=M$). Moreover, proposition \ref{main_prop} resp.~2 shows that the map $f:(U,g)\to (f(U),g_0)$ is a conformal diffeomorphism where the conformal factor is given by the function $\alpha =\sqrt{\mathrm{area}(f)}$, i.e.~$(U,\overline{g}:=\alpha g)$ and $(f(U),g)$ are locally isometric. We will use equation (\ref{scal_equal}) to prove that $\alpha $ has to be constant which implies that $U=M$ and that $f:(M,\alpha g)\to (M_0,g_0)$ must be a Riemannian covering (note that $f$ is surjective since $\deg (f)\neq 0$). We first note that equation (\ref{scal_equal}) together with $\mathrm{scal}_0\circ f>0$ imply that $\alpha $ is smooth on all of $M$ not only on $U$. Moreover, the scalar curvature of $g$ and $\overline{g}=\alpha g$ are on $U$ related by
\[
\overline{\mathrm{scal}}=\frac{1}{\alpha }\mathrm{scal}_g+\frac{n-1}{\alpha ^2}\delta \mathrm{d}\alpha -\frac{(n-1)(n-6)}{4\alpha ^3}|\mathrm{d}\alpha |^2_g.
\]
Since $\overline{g}=f^*g_0$ on $U$, the scalar curvature of $\overline{g}$ is given by $\overline{\mathrm{scal}}=\mathrm{scal}_0\circ f$ and we conclude from equation (\ref{scal_equal})
\[
\frac{n-1}{\alpha ^2}\delta \mathrm{d}\alpha -\frac{(n-1)(n-6)}{4\alpha ^3}|\mathrm{d}\alpha |^2_g=0.
\]
Thus, since $\alpha $ is smooth on $M$ [cf.~equation (\ref{scal_equal})] and $\alpha =0$ on $M-U$, the following equation holds on all of $M$ for all $k\geq 1$:
\[
\alpha ^{k}\delta \mathrm{d}\alpha -\frac{n-6}{4}\alpha ^{k-1}|\mathrm{d}\alpha |^2=0.
\]
Integrate over $M$ w.r.t.~the volume form of $g$ yields for all $k\geq 1$
\[
\left( k-\frac{n-6}{4}\right) \int\limits _M\alpha ^{k-1}|\mathrm{d}\alpha |^2=0
\]
and hence, $\mathrm{d}\alpha =0$ shows $\alpha =\sqrt{\mathrm{area}(f)}=\mathrm{const}$ which completes the proof of the main theorems (recall we have $\alpha >0$ since $\deg (f)\neq 0$ yields $U\neq \emptyset $). Corollary 1 follows simply by the fact that a locally symmetric space with Ricci curvature $\mathrm{Ric}_0>0$ has non--negative curvature operator.
\section{Conformal submersions and scalar curvature}

Let $f:M\to M_0$ be a smooth map between closed oriented manifolds. The \emph{$\widehat{A}$--degree} of $f$ is defined as
\[
\deg _{\widehat{A}}(f):=\int \limits _M \widehat{\mathbf{A}}(M)\cdot f^*\omega 
\]
where $\widehat{\mathbf{A}}(M)$ denotes the total $\widehat{A}$--class of $M$ and $\omega $ is a volume form on $M_0$ with $\int _{M_0}\omega =1$. In particular, if $\deg _{\widehat{A}}(f)$ is non--zero, then $\dim M= \dim M_0+4k$ with $k\geq 0$ and if $\dim M=\dim M_0$, then $\deg _{\widehat{A}}(f)=\deg (f)$. The $\widehat{A}$--degree of $f$ can also be defined in terms of regular values of $f$. If $\dim M-\dim M_0=4k$, $\mathrm{deg}_{\widehat{A}}(f)$ equals the $\widehat{A}$--genus of the manifold $f^{-1}(y)$ where $y\in M_0$ is a regular value. A smooth map $f:(M,g)\to (M_0,g_0)$ is called a \emph{conformal submersion} if there is a smooth function $\alpha :M\to (0,\infty )$ such that $f:(M,\overline{g}:=\alpha g)\to (M_0,g_0)$ is a Riemannian submersion (cf.~\cite{Bes}). The function $\alpha $ is said to be the \emph{conformal factor} of the conformal submersion. We denote by $A$ and $T$ the invariants introduced in \cite[Ch.~9]{Bes} for the Riemannian submersion $f:(M,\overline{g})\to (M_0,g_0)$. Moreover let $N$ be the mean curvature vector field of the submanifolds $f^{-1}(y)\subset M$, then $N$ is given by the $\overline{g}$--trace of the second fundamental form $T$. 
\begin{thm}
Let $(M_0^m,g_0)$ be an oriented closed Riemannian manifold of dimension $m\geq 3$ and with $\curv{R}^0\geq 0$, $\mathrm{Ric}_0>0$ and $\chi (M_0)\neq 0$. Suppose $(M^n,g)$ is an oriented closed manifold with $\mathrm{Ric}_g(p)\not\equiv 0$ for all $p\in M$ and $f:M\to M_0$ is a spin map of nonzero $\widehat{A}$--degree: $\deg _{\widehat{A}}(f)\neq 0$. If the scalar curvature satisfies
\begin{equation}
\label{main_in}
\mathrm{scal}_g\geq \sqrt{\mathrm{area}(f)}\cdot \bigl( \mathrm{scal}_0\circ f\bigl) ,
\end{equation}
then inequality (\ref{main_in}) is an equality and $f:(M,g)\to (M_0,g_0)$ is a conformal submersion where the conformal factor is given by $\alpha :=\sqrt{\mathrm{area}(f)}$. If the mean curvature vector field $N$ satisfies $|N|^2_{\overline{g}}\leq |T|^2_{\overline{g}}$, then the conformal factor $\alpha $ is constant and the invariant $A$ vanishes (i.e.~the horizontal space $\bundle{H}=\ker (f_*)^\perp \subset TM$ is integrable). Furthermore, if the fibers of this conformal submersion are all minimal, i.e.~$N=0$, then the conformal factor $\alpha $ is constant and $(M,\overline{g}=\alpha g)$ is locally isometric to $(M_0\times X,g_0\oplus h)$ where $(X,h)$ is a closed Ricci flat spin manifold of nonzero $\widehat{A}$--genus. 
\end{thm}

The assumption $\mathrm{Ric}_g(p)\not \equiv 0$ can be replaced by $\mathrm{area}(f)>0$. The proof of this version of the theorem needs $\mathrm{area}(f)>0$ in order to conclude that $f:M\to M_0$ is a submersion. If $\dim M=\dim M_0$ we used a simple conformal argument to show that $\mathrm{area}(f)$ has to be a nonzero constant.

\begin{proof}
Let $\bundle{V}^{2k}$ be an oriented Riemannian vector bundle over an oriented Riemannian manifold $(M^{2n},g)$ in such a way that $\bundle{E}=TM\oplus \bundle{V}$ admits a spin structure. Then the complex spinor bundle $\spinor \bundle{E}$ is naturally $\Z _2$ graded by the volume form of $\Cl _\C (\bundle{E})$. Moreover, using the embedding $\Cl _\C (TM)\hookrightarrow \Cl _\C (\bundle{E})$ and an arbitrary Riemannian connection on $\bundle{V}$, $\spinor \bundle{E}$ becomes a complex Dirac bundle over $M$ and the index of the Dirac operator $\overline{\dirac }^+:\Gamma (\spinor ^+\bundle{E})\to \Gamma (\spinor ^-\bundle{E})$ is given by
\[
\mathrm{ind}(\overline{\dirac }^+)=\{ \boldsymbol \chi (\bundle{V})\widehat{\mathbf{A}}(\bundle{V})^{-1}\widehat{\mathbf{A}}(M)\} [M]
\]
where $\boldsymbol \chi (\bundle{V})\in H^{2k}(M;\Q )$ is the Euler class and $\widehat{\mathbf{A}}(\bundle{V})\in H^{4*}(M;\Q )$ is the $\widehat{A}$--class of $\bundle{V}$. This version of the Atiyah Singer index theorem can be deduced from Chapter III theorem 13.13 and proposition 11.24 in \cite{LaMi}. 

In the situation of theorem 3, $\chi (M_0)\neq 0$ yields $\dim M_0$ is even and $\deg _{\widehat{A}}(f)\neq 0$ yields $\dim M=\dim M_0+4k$. In particular, we obtain in case $\bundle{V}=f^*TM_0$
\[
\begin{split}
\boldsymbol \chi (\bundle{V})\widehat{\mathbf{A}}(\bundle{V})^{-1}=&f^*(\boldsymbol \chi (TM_0)\widehat{\mathbf{A}}(TM_0)^{-1})\\
=& f^*(\boldsymbol \chi (\bundle{V}))=\chi (M_0)f^*\omega
\end{split} 
\]
where $\omega $ means the orientation form. Hence, the index of the Dirac operator $\overline{\dirac }^+:\Gamma (\spinor ^+\bundle{E})\to \Gamma (\spinor ^-\bundle{E})$ is determined by
\[
\mathrm{ind}(\overline{\dirac }^+)=\chi (M_0)\{ \widehat{\mathbf{A}}(M)f^*\omega \} [M]= \chi (M_0)\cdot \deg _{\widehat{A}}(f) \neq 0
\]
and $\ker (\overline{\dirac})$ is non--trivial. Thus, the integrated Bochner--Weitzenb\"ock formula (\ref{int_lich}) together with inequalities (\ref{ineq}) and (\ref{main_in}) prove equality in (\ref{main_in}) and the fact that the minimal eigenvalue of $\frak{R}$ is at each point given by $-\frac{1}{4}\sqrt{\mathrm{area}(f)}\cdot \mathrm{scal}_0\circ f$. As usual set $\alpha :=\sqrt{\mathrm{area}(f)}$ and consider the open set $U=\{ p\in M|\alpha (p)>0\} $, then proposition \ref{main_prop} shows that $f:(U,\alpha g)\to (f(U),g_0)$ is a Riemannian submersion (note that $U$ is not empty because $\deg _{\widehat{A}}(f)\neq 0$). The next step will show that $U=M$. Let $0\neq \psi \in \ker (\overline{\dirac})$, then $\psi $ is parallel w.r.t.~$\overline{\nabla }$ and $\overline{R}_{x,y}\psi =0$ yields:
\[
0=(R_{x,y}\otimes \mathrm{Id} +\mathrm{Id}\otimes R^0_{f_*x,f_*y})\psi .
\]
Hence, Clifford multiplication supplies
\[
\mathrm{ric}_g(y)\cdot \psi =2\sum _{j=1}^n(e_j\cdot R_{e_j,y}\otimes \mathrm{Id})\psi =-2\sum _{j=1}^n(e_j\otimes R^0_{f_*e_j,f_*y})\psi .
\]
If $\alpha $ vanishes at $p\in M$, so does $f_*:\Lambda ^2T_pM\to \Lambda ^2T_{f(p)}M_0$ which proves in particular that $\mathrm{ric}_g(y)\cdot \psi =0$ and thus, $\mathrm{ric}_g(y)=0 $ for all $y\in T_pM$, i.e.~we obtain $\mathrm{Ric}_g(p)=0$ for all $p\in M-U$. However, we assumed $\mathrm{Ric}_g(p)\neq 0$ for all $p\in M$, i.e.~$U=M$ and $f:(M,g)\to (M_0,g_0)$ is a conformal submersion while the conformal factor is determined by $\alpha $ (note that $f$ is surjective since $\deg _{\widehat{A}}(f)\neq 0$, and $\alpha $ is smooth since inequality (\ref{main_in}) is an equality and $\mathrm{scal}_0\circ f>0$). This completes the first part of theorem 3.

Consider the orthogonal decomposition $TM=\ker (f_*)\oplus \bundle{H}$ over $M$, then $\ker (f_*)$ is an integrable distribution and $\bundle{H}$ is integrable if and only if the invariant $A$ vanishes. If $y\in \ker (f_*)$ at an arbitrary point of $M$, the right hand side of the above equation for the Ricci tensor vanishes, and so $\mathrm{ric}_g(y)=0$. Hence, we conclude that the Ricci tensor of $g$ restricted to $\ker (f_*)$ is trivial. Let $\overline{g}:=\alpha g$ be the conformal transformation of $g$ on $M$ by $\alpha =\sqrt{\mathrm{area}(f)}$ (we already proved that $\alpha >0$ on $M$). Then the Ricci tensors of $\overline{g}$ and $g$ are related by (use derivatives w.r.t.~$g$)
\[
\overline{\mathrm{Ric}}=\mathrm{Ric}-\frac{n-2}{2\alpha }\nabla ^2\alpha +\frac{3(n-2)}{4\alpha ^2}\mathrm{d}\alpha \otimes \mathrm{d}\alpha +\left[ \frac{1}{2\alpha }\delta \mathrm{d}\alpha -\frac{n-4}{4\alpha ^2}|\mathrm{d}\alpha |^2_g\right] g
\]
and the scalar curvature satisfies
\[
\overline{\mathrm{scal}}=\frac{1}{\alpha }\mathrm{scal}_g+\frac{n-1}{\alpha ^2}\delta \mathrm{d}\alpha -\frac{(n-1)(n-6)}{4\alpha ^3}|\mathrm{d}\alpha |^2_g.
\]
Taking the $\overline{g}$--trace of $\overline{\mathrm{Ric}}_{|\bundle{H}}$ (restriction of $\overline{\mathrm{Ric}}$ to the horizontal distribution $\bundle{H}$) leads to
\[
\begin{split}
\mathrm{tr}_{\overline{g}}(\overline{\mathrm{Ric}}_{|\bundle{H}})&=\sum _{j=1}^m\overline{\mathrm{Ric}}(x_j,x_j)\\
&=\frac{1}{\alpha }\mathrm{scal}_g+\frac{m}{2\alpha ^2}\delta \mathrm{d}\alpha -\frac{m(n-4)}{4\alpha ^3}|\mathrm{d}\alpha |^2_g+\frac{3(n-2)}{4\alpha ^3}|\pi _{\bundle{H}}\mathrm{d}\alpha | _g^2-\frac{n-2}{2\alpha }F(\alpha )
\end{split}
\]
where
\[
F(\alpha )=\sum _{j=1}^k\nabla ^2\alpha (x_j,x_j) 
\]
for an $\overline{g}$--orthonormal basis $x_1,\ldots ,x_m$ of $\bundle{H}$. In this case we used the fact that $\frac{x_j}{\sqrt{\alpha }}$, $j=1\ldots m$, is a $g$--orthonormal basis of $\bundle{H}$, $\mathrm{rk}(\bundle{H})=\dim M_0=m$ and that the Ricci tensor of $g$ is trivial on $\bundle{H}^\perp =\ker (f_*)$ as shown above. Since $f:(M,\overline{g})\to (M_0,g_0)$ is a Riemannian submersion, we obtain a relation of the Ricci tensors of $\overline{g}$ and $g_0$ (cf.~equation (9.36c) in \cite{Bes}, we use the same notation in order to avoid any confusion, in particular $D_X$ is the Levi--Civita connection of $\overline{g}$):
\[
\overline{\mathrm{Ric}}(X,Y)=\mathrm{Ric}_0(f_*X,f_*Y)-2(A_X,A_Y)-(TX,TY)+\frac{1}{2}((D_XN,Y)+(D_YN,X))
\]
where $X,Y\in \bundle{H}$. Taking the $\overline{g}$--trace of this formula yields
\[
\mathrm{tr}_{\overline{g}}(\overline{\mathrm{Ric}}_{|\bundle{H}})=\mathrm{scal}_0\circ f-2|A|^2_{\overline{g}}-|T|^2_{\overline{g}}-\breve{\delta }N.
\]
Combining the two equations for $\mathrm{tr}_{\overline{g}}(\overline{\mathrm{Ric}}_{|\bundle{H}})$ and using $\mathrm{scal}_g=\alpha \cdot \mathrm{scal}_0\circ f$, we obtain
\[
\frac{m}{2\alpha ^2}\delta \mathrm{d}\alpha -\frac{m(n-4)}{4\alpha ^3}|\mathrm{d}\alpha |^2_g+\frac{3(n-2)}{4\alpha ^3}|\pi _{\bundle{H}}\mathrm{d}\alpha | _g^2-\frac{n-2}{2\alpha }F(\alpha )=-2|A|^2_{\overline{g}}-|T|^2_{\overline{g}}-\breve{\delta }N.
\]
The next step is to simplify $F(\alpha )$ and $\breve{\delta }N$. Let $\overline{\delta }$, $\delta $ be the (negative) divergence operator of $\overline{g}$, $g$ respectively, then $\overline{\delta }V=\delta V-\frac{n}{2\alpha }\mathrm{d}\alpha (V)$ for arbitrary vector fields $V$. Moreover, if $X_1,\ldots ,X_m$ is an $\overline{g}$--orthonormal frame of $\bundle{H}$ and $U_1,\ldots ,U_{n-m}$ is an $\overline{g}$--orthonormal frame of $\ker (f_*)$, we obtain from $N\in \Gamma (\bundle{H})$
\[
\overline{\delta }N=-\sum _{j=1}^m\overline{g}(D_{X_j}N,X_j)-\sum _{i=1}^{n-m}\overline{g}(D_{U_j}N,U_j)=\breve{\delta }N+\sum _{i=1}^{n-m}\overline{g}(N,D_{U_i}U_i)
\]
In particular, the definition of $T$ and $N$ in \cite[Chapter 9]{Bes} supply
\[
\breve{\delta } N=\delta N-\frac{n}{2\alpha }\mathrm{d}\alpha (N)-|N|^2_{\overline{g}}.
\]
Moreover, using local frames for $\bundle{H}$ and $\ker (f_*)$ a straightforward calculation shows
\[
F(\alpha )=\mathrm{tr}_{\overline{g}}(\nabla ^2\alpha _{|\bundle{H}})=-\frac{1}{\alpha }\delta \pi _{\bundle{H}}\mathrm{d}\alpha + \mathrm{d}\alpha (N)+\frac{n-m}{2\alpha ^2}| \pi _{\bundle{H}}\mathrm{d}\alpha |^2_g-\frac{m}{2\alpha ^2}|\pi _{\bundle{H}^\perp }\mathrm{d}\alpha |^2_g
\]
where $N=\sum T_{U_i}U_i\in \Gamma (\bundle{H})$ (cf.~\cite[9.34]{Bes}) and we used that $\sum A_{X_i}X_i=0$ for a $\overline{g}$--orthonormal frame $X_1,\ldots ,X_m\in \Gamma (\bundle{H})$. Insert the expressions for $F(\alpha )$ and $\breve{\delta }N$ into the above equation yields
\begin{equation}
\label{last_eq}
\begin{split}
\frac{m}{2\alpha ^2}\delta \mathrm{d}\alpha +\frac{n-2}{2\alpha ^2}\delta \pi _\bundle{H}\mathrm{d}\alpha &- \frac{n^2-5n-2m+6}{4\alpha ^3}|\pi _{\bundle{H}}\mathrm{d}\alpha | _g^2+\frac{2m}{4\alpha ^3}|\pi _{\bundle{H}^\perp }\mathrm{d}\alpha |^3_g=\\
&=-2|A|^2_{\overline{g}}-|T|^2_{\overline{g}}+|N|^2_{\overline{g}}-\delta N+\frac{n-1}{\alpha }\mathrm{d}\alpha (N).
\end{split}
\end{equation}
Multiply this equation by $\alpha ^{n-1}$ and integrate with respect to $\mathrm{vol}_g$, then the last two terms on the right hand side of this equation vanish, since $\int \alpha ^{n-1}\delta N=(n-1)\int \alpha ^{n-2}\mathrm{d}\alpha (N)$, in fact we obtain
\[
\int \limits _Mc(n,m)\cdot \alpha ^{n-4}|\pi _\bundle{H}\mathrm{d}\alpha |^2_g+\frac{m(n-2)}{2}\alpha ^{n-4}|\pi _{\bundle{H}^\perp }\mathrm{d}\alpha |^2_g=-\int \limits _M 2|A|^2 _{\overline{g}}+\left[ |T|^2 _{\overline{g}}-|N|^2_{\overline{g}}\right] 
\]
where $c(n,m)$ is defined by
\[
c(n,m):=\frac{n^2-5n+6+2m(n-2)}{4}.
\]
In case $n,m\geq 3$, $c(n,m)$ is clearly positive. Thus, the left hand side of the above equation is non--negative while the right hand side is non--positive under the assumption $|N|^2_{\overline{g}}\leq |T|^2_{\overline{g}}$ which means that $\mathrm{d}\alpha =0$, $|A|^2_{\overline{g}}=0$ and $|T|^2_{\overline{g}}=|N|^2_{\overline{g}}$. The remaining statements are conclusion from \cite[Ch.~9]{Bes}. Note that if the fibers are all minimal, we also conclude $|T|^2_{\overline{g}}=0$, i.e.~$(M,g)$ is locally a Riemannian product. The fact that $(X,h)$ is Ricci flat follows from $\overline{\mathrm{Ric}}=\mathrm{Ric}_g$ and the fact that $\mathrm{Ric}_g$ is trivial on $\ker (f_*)$. The $\widehat{A}$--genus of $X$ is nonzero since $\widehat{A}(X)=\deg _{\widehat{A}}(f)$. 
\end{proof}
\bibliographystyle{abbrv}
\bibliography{bibliothek}

\def\cprime{$'$}
\begin{thebibliography}{10}

\bibitem{Bes}
A.~L. Besse.
\newblock {\em Einstein manifolds}, volume~10 of {\em Ergebnisse der Mathematik
  und ihrer Grenzgebiete (3) [Results in Mathematics and Related Areas (3)]}.
\newblock Springer-Verlag, Berlin, 1987.

\bibitem{Go1}
S.~Goette.
\newblock Scalar curvature estimates by parallel alternating torsion.
\newblock {\em arXiv: 0709.4586[math.DG]}, 2007.

\bibitem{GoSe1}
S.~Goette and U.~Semmelmann.
\newblock Scalar curvature estimates for compact symmetric spaces.
\newblock {\em Differential Geom. Appl.}, 16(1):65--78, 2002.

\bibitem{Gr01}
M.~Gromov.
\newblock Positive curvature, macroscopic dimension, spectral gaps and higher
  signatures.
\newblock In {\em Functional analysis on the eve of the 21st century, Vol.\ II
  (New Brunswick, NJ, 1993)}, volume 132 of {\em Progr. Math.}, pages 1--213.
  Birkh\"auser Boston, Boston, MA, 1996.

\bibitem{GrLa2}
M.~Gromov and H.~B. Lawson, Jr.
\newblock Spin and scalar curvature in the presence of a fundamental group.
  {I}.
\newblock {\em Ann. of Math. (2)}, 111(2):209--230, 1980.

\bibitem{Kra}
W.~Kramer.
\newblock The scalar curvature on totally geodesic fiberings.
\newblock {\em Ann. Global Anal. Geom.}, 18(6):589--600, 2000.

\bibitem{LaMi}
H.~B. Lawson, Jr. and M.-L. Michelsohn.
\newblock {\em Spin geometry}, volume~38 of {\em Princeton Mathematical
  Series}.
\newblock Princeton University Press, Princeton, NJ, 1989.

\bibitem{Llarull}
M.~Llarull.
\newblock Sharp estimates and the {D}irac operator.
\newblock {\em Math. Ann.}, 310(1):55--71, 1998.

\bibitem{Loh2}
J.~Lohkamp.
\newblock Metrics of negative {R}icci curvature.
\newblock {\em Ann. of Math. (2)}, 140(3):655--683, 1994.

\bibitem{MinO2}
M.~Min-Oo.
\newblock Scalar curvature rigidity of certain symmetric spaces.
\newblock In {\em Geometry, topology, and dynamics (Montreal, PQ, 1995)},
  volume~15 of {\em CRM Proc. Lecture Notes}, pages 127--136. Amer. Math. Soc.,
  Providence, RI, 1998.

\end{thebibliography}

\end{document}